\documentclass[graybox]{svmult}

\usepackage{mathptmx}         
\usepackage{helvet}           
\usepackage{courier}          
\usepackage{type1cm}          
\usepackage{makeidx}          
\usepackage{graphicx}         
\usepackage{multicol}         
\usepackage[bottom]{footmisc} 
\usepackage[T1]{fontenc}      
\usepackage{relsize}          

\makeindex             

\newif\iffinal

\finaltrue

\usepackage[caption=false]{subfig}

\usepackage{dsfont}

\usepackage{amsmath}
\usepackage{amssymb}

\DeclareMathOperator*{\argmax}{arg\,max}
\DeclareMathOperator {\dist}  {dist}
\DeclareMathOperator {\im}    {im}

\DeclareMathOperator {\rest}  {\phi} 

\newcommand{\card}[1]         {\ensuremath{|#1|}}

\newcommand{\Li}              {\ensuremath{\mathrm{L}_\infty}}
\newcommand{\manifold}        {\ensuremath{\mathds{M}}}
\newcommand{\real}            {\ensuremath{\mathds{R}}}

\newcommand{\betti}                  [1]{\ensuremath{\beta_{#1}}}
\newcommand{\boundary}               [1]{\ensuremath{\partial_{#1}}}
\newcommand{\boundarygroup}          [1]{\ensuremath{B_{#1}}}
\newcommand{\chaingroup}             [1]{\ensuremath{C_{#1}}}
\newcommand{\cyclegroup}             [1]{\ensuremath{Z_{#1}}}
\newcommand{\homologygroup}          [1]{\ensuremath{H_{#1}}}
\newcommand{\persistenthomologygroup}[2]{\ensuremath{\homologygroup{#1}^{#2}}}
\newcommand{\simplicialcomplex}         {\mathrm{K}}

\usepackage{tikz}
\usepackage{pgfplots}

\pgfplotsset{compat=1.13}

\usepgfplotslibrary{colorbrewer}
\usepgfplotslibrary{external}
\usepgfplotslibrary{fillbetween}

\tikzsetexternalprefix{Figures/External/}

\definecolor{amber}     {RGB}{255,191,  0}
\definecolor{cardinal}  {RGB}{196, 30, 58}
\definecolor{yale}      {RGB}{ 70,130,180}

\pgfplotscreateplotcyclelist{color list}{%
  cardinal,
  yale,
  amber
}

\pgfplotsset{cycle list name=color list}

\usetikzlibrary{calc}
\usetikzlibrary{shapes}

\makeatletter
\newcommand{\twodigits}[1]{\two@digits{#1}}
\makeatother

\usepackage[final]{changes}
\definechangesauthor[name={BR}, color=blue]{br}

\usepackage[olditem,
            oldenum]{paralist}

\usepackage{xspace}

\renewcommand{\th}{\textsuperscript{\textup{th}}\xspace}

\usepackage[absolute]{textpos}

\begin{document}

\title*{Persistent Intersection Homology for the Analysis of Discrete Data}

\author{Bastian Rieck \and Markus Banagl \and Filip Sadlo \and Heike Leitte}

\institute{Bastian Rieck \and Heike Leitte \at TU Kaiserslautern, \email{\{rieck, leitte\}@cs.uni-kl.de}
  \and Markus Banagl \and Filip Sadlo \at Heidelberg University, \email{banagl@mathi.uni-heidelberg.de}, \email{sadlo@uni-heidelberg.de}%
}

\iffinal
\else
  \begin{textblock*}{\paperwidth}[1, 0](\paperwidth,0.5cm)
  \scriptsize%
  Authors' copy. Please refer to \emph{Topological Methods in Data Analysis and Visualization
  V: Theory, Algorithms, and Applications} for the definitive version of this chapter.
  \end{textblock*}
\fi

\maketitle

\abstract{%
  Topological data analysis is becoming increasingly relevant to support the analysis of
  unstructured data sets. A common assumption in data analysis is that the data set is
  a sample---not necessarily a uniform one---of some high-dimensional manifold.
  In such cases, persistent homology can be successfully employed to extract features, remove noise,
  and compare data sets. The underlying problems in some application domains, however, turn out to
  represent \emph{multiple} manifolds with different dimensions.
  Algebraic topology typically analyzes such problems using intersection homology, an extension of
  homology that is capable of handling configurations with singularities. In this paper, we describe
  how the persistent variant of intersection homology can be used to assist data analysis in
  visualization. We point out potential pitfalls in approximating data sets with singularities
  and give strategies for resolving them.
}

\section{Introduction}

The \emph{manifold hypothesis} is a traditional assumption for the analysis of multivariate data.
Briefly put, it assumes that the input data are a sample of some manifold~\manifold, whose
intrinsic dimension~$d$ is much smaller than the ambient dimension~$D$.
Typical examples of this assumption are found in dimensionality reduction algorithms~\cite{Roweis00,
Tenenbaum00}.
For certain applications, such as image analysis~\cite{Donoho05} or image
recognition~\cite{Hinton97}, we already know this hypothesis to be true---at least with respect to
the models that are often used to describe such data.
For other applications, there are strategies~\cite{Fefferman16, Narayanan10} for testing this
hypothesis provided that a sufficient number of samples is available.

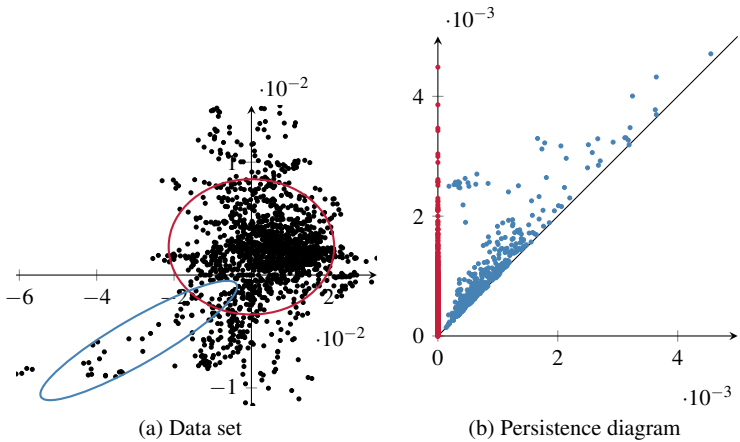
\begin{figure}[tbp]
  \center
  \subfloat[Data set\label{sfig:Data set example}]{%
    \iffinal
      \includegraphics{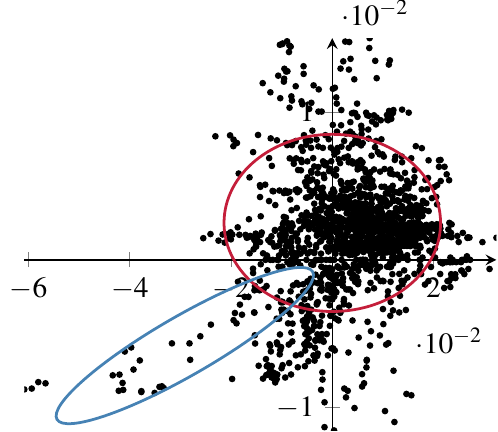}
    \else
      \begin{tikzpicture}
        \begin{axis}[%
          axis lines        = middle,
          mark size         = 0.75pt,
          scale             = 0.70,
          clip marker paths = true,
        ]
          \addplot[only marks] file {Data/MMM/SAP_CPI_02.txt};
          \draw[thick, cardinal      ] ( 0.000, 0.0025) ellipse (1.10cm and 0.9cm);
          \draw[rotate=30,thick, yale] (-0.025,-0.0120) ellipse (1.50cm and 0.3cm);
        \end{axis}
      \end{tikzpicture}
    \fi
  }
  \quad
  \subfloat[Persistence diagram\label{sfig:Persistence diagram example}]{%
    \iffinal
      \includegraphics{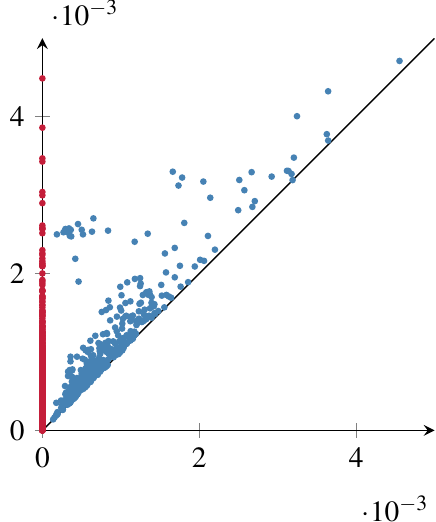}
    \else
      \begin{tikzpicture}
        \begin{axis}[%
          axis x line       = bottom,
          axis y line       = left,
          enlargelimits     = false,
          mark size         = 0.75pt,
          scale             = 0.70,
          unit vector ratio*= 1 1 1,
        ]
          \addplot gnuplot[only marks, raw gnuplot] {%
            plot "Data/MMM/SAP_CPI_02_persistence_diagrams.txt" index 0 with points
          };
          \addplot gnuplot[only marks, raw gnuplot] {%
            plot "Data/MMM/SAP_CPI_02_persistence_diagrams.txt" index 1 with points
          };

          \addplot[domain=0.0:0.005] {x};
        \end{axis}
      \end{tikzpicture}
    \fi
  }
  \caption{
    \protect\subref{sfig:Data set example} The structure of a central ``core'' with ``flares''
    emanating from it appears in many data sets~(here, 2-year growth rates of Standard~\&~Poor's 500
    vs.\ the U.S.\ CPI with the core shown in red and one example flare shown in blue).
    \protect\subref{sfig:Persistence diagram example} The corresponding persistence diagram shows
    topological features in dimension zero~(red) and dimension one~(blue).
  }
  \label{fig:Market data}
\end{figure}

The \emph{practice} of multivariate data analysis seems to suggest something else, though:
Carlsson~\cite{Carlsson14}, for example, remarks that many real-world data sets exhibit a central
``core'' structure, from which different ``flares'' emanate.
Figure~\ref{fig:Market data} illustrates this for a simple 2D data set, generated from 2-year growth
rates of Standard~\&~Poor's 500 vs.\ the U.S.\ CPI.
This structure is irreconcilable with the structure of a single manifold.
Novel data analysis algorithms such as \textsc{Mapper}~\cite{Singh07} account for this fact by not
making any assumptions about manifold structures and attempting to fit data in a local manner---a
strategy that is also employed in low-dimensional manifold learning~\cite{Saul03}.

In this paper, we argue that some real-world data sets require special tools to assess their
structure. Just as persistent homology~\cite{Edelsbrunner10,Edelsbrunner14} was originally developed
to analyze samples from spaces that are supposed to have the structure of a manifold, we need
a special tool to analyze spaces for which this assumption does \emph{not} hold.
More precisely, we will tackle the task of analyzing spaces that are composed of different
manifolds~(with possibly varying dimensions) using \emph{intersection homology}~\cite{Kirwan06} and
\emph{persistent intersection homology}~\cite{Bendich09,Bendich11}.
To make it accessible to a wider community of researchers, we devote a large portion of this paper
to explaining the theory behind persistent intersection homology.
Furthermore, we discuss implementation details and present an open-source framework for its
calculation. We also describe pitfalls in ``naive'' applications of persistent intersection homology
and develop strategies to resolve them.

\section{Background}

We first explain the mathematical tools required to describe spaces that are not composed of
a single manifold, but of multiple ones. Next, we introduce~(persistent) intersection homology, give
a brief algorithm for its computation, and describe how to use it to analyze real-world data
sets.

\subsection{Stratifications}
\label{sec:Stratifications}

\begin{figure}[tbp]
  \centering
  \subfloat[][\label{fig:Pinched torus}]{%
    \includegraphics[width=4cm]{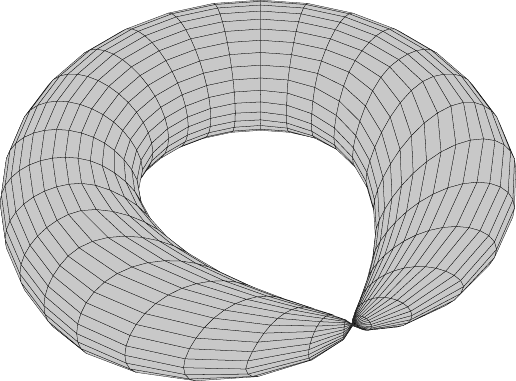}
  }
  \quad
  \subfloat[][\label{fig:Pinched torus curvature}]{%
    \includegraphics[width=4cm]{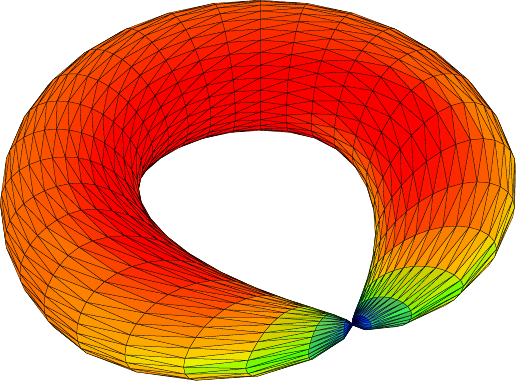}
  }
  \caption{%
    \protect\subref{fig:Pinched torus} The ``pinched torus'' is a classical example of an object
    that is \emph{not} a manifold but composed of parts that are manifolds, provided the singular
    point that is caused by the ``pinch'' is ignored.
    \protect\subref{fig:Pinched torus curvature} The singular point is readily visible when
    calculating mean curvature estimates.
  }
\end{figure}

Stratifications are a way of describing spaces that are not a manifold per se, but composed of
multiple parts, each of which is a manifold. A common example of such a space is the ``pinched
torus'', which is obtained by collapsing~(i.e., pinching) one minor ring of the torus to a single
point. Figure~\ref{fig:Pinched torus} depicts an example. The neighborhood of the pinch point is
\emph{singular} because it does not satisfy the conditions of a manifold: it does not have
a neighborhood that is homeomorphic to a ball. If we remove this singular point, however, the
remaining space is just a~(deformed) cylinder, i.e., a manifold.
Permitting the removal of certain parts of a space may thus be beneficial to describe the
manifolds it is composed of. This intuition leads to the concept of stratifications.

Let $X\subseteq\real^n$ be a topological space. A \emph{topological stratification} of $X$ is a filtration of closed subspaces
\begin{equation}
  \emptyset \subseteq X_{-1} \subseteq X_0 \subseteq X_1 \subseteq \dots \subseteq X_{d-1} \subseteq X_d = X,
\end{equation}
such that for each $i$ and every point $x \in X_i \setminus X_{i-1}$ there is a neighborhood $U
\subseteq X$ of $x$, a compact $(n-1-i)$-dimensional stratified topological space V, and
a filtration-preserving homeomorphism $U \simeq \real^i \times CV$, where $CV$ denotes the
\emph{open cone} on $V$, i.e., $CV := V \times [0,1) / ( V \times \{0\} )$.
We refer to $X_i \setminus X_{i-1}$ as the $i$-dimensional \emph{stratum} of $X$. Notice that it is
always a~(smooth) manifold, even though the original space might not be a manifold.
Hence, this rather abstract definition turns out to be a powerful description for a large family of
spaces.
There are some stratifications with special properties that are particularly suited for analyzing
spaces. Goresky and MacPherson~\cite{Goresky80}, the inventors of intersection homology,
suggest using a stratification that satisfies $X_{d-1} = X_{d-2}$ so that the
$(d-1)$-dimensional stratum is empty, i.e., $X_{d-1} \setminus X_{d-2} = \emptyset$.

\subsection{Homology and Persistent Homology}

Prior to introducing~(persistent) intersection homology, we briefly describe simplicial homology and
its persistent counterpart. Given a $d$-dimensional simplicial complex~$\simplicialcomplex$, the
chain groups $\{\chaingroup{0}$, \dots, $\chaingroup{d}\}$ contain formal sums~(simplicial chains)
of simplices of a given dimension.
A boundary operator $\boundary{p}\colon\chaingroup{p}\to\chaingroup{p-1}$ satisfying
$\boundary{p-1}\circ\boundary{p} = 0$~(i.e., a closed boundary does not have a boundary itself) then
permits us to create a chain complex from the chain groups.
This results in two subgroups, namely the cycle group~$\cyclegroup{p} := \ker\boundary{p}$ and the
boundary group~$\boundarygroup{p} := \im\boundary{p+1}$,
from which we obtain the $p$\th homology group as
\begin{equation}
    \homologygroup{p} := \cyclegroup{p} / \boundarygroup{p},
\end{equation}
where the $/$-operator refers to the quotient group. Intuitively, elements in the cycle
group~$\cyclegroup{p}$ constitute sets of simplicial chains that do not have a boundary, while
elements in the boundary group~$\boundarygroup{p}$ are the boundaries of higher-dimensional
simplices. By removing these in the definition of the homology group, we obtain a group that
describes high-dimensional ``holes'' in~$\simplicialcomplex$.

Homology is a powerful tool to discriminate between different triangulated topological spaces. It is
common practice to use the Betti numbers~$\betti{p}$, i.e., the ranks of the homology groups, to
obtain a signature of a space. In practice, the Betti numbers turn out to be highly susceptible to
noise, which prompted the development of persistent homology~\cite{Edelsbrunner10}.
Its basic premise is that the simplicial complex~$\simplicialcomplex$ is associated with
a filtration,
 \begin{equation}
  \emptyset = \simplicialcomplex_0 \subseteq \simplicialcomplex_1 \subseteq \dots \subseteq \simplicialcomplex_{n-1} \subseteq \simplicialcomplex_n = \simplicialcomplex,
\end{equation}
where each $\simplicialcomplex_i$ is typically assigned a function value, such as a distance. The
filtration induces a homomorphism of the corresponding homology groups, i.e.,
$f_p^{i,j} \colon \homologygroup{p}(\simplicialcomplex_i) \to \homologygroup{p}(\simplicialcomplex_j)$,
leading to the definition of the $p$\th persistent homology group~$\persistenthomologygroup{p}{i,j}$
for two indices $i \leq j$ as
\begin{equation}
  \persistenthomologygroup{p}{i,j} := \cyclegroup{p}\left(\simplicialcomplex_i\right) / \left( \boundarygroup{p}\left(\simplicialcomplex_j\right)\cap\cyclegroup{p}\left(\simplicialcomplex_i\right)\right).
\end{equation}
This group contains all the homology classes of~$\simplicialcomplex_i$ that are still
present in~$\simplicialcomplex_j$. It is possible to keep track of all homology classes within the
filtration.

The calculation of persistent homology results in a set of pairs~$(i,j)$, which denote
a homology class that was created in~$\simplicialcomplex_i$ and destroyed~(vanished) in~$\simplicialcomplex_j$.
Letting $f_i$ denote the associated function value of $\simplicialcomplex_i$, these pairs are
commonly visualized in a persistence diagram~\cite{Cohen-Steiner07} as~$(f_i, f_j)$. The distance of
each pair to the diagonal, measured in the $\Li$-norm, is referred to as the \emph{persistence} of
a topological feature.
It is now common practice in topological data analysis to use persistence to separate noise from
salient features in real-world data sets~\cite{Edelsbrunner10,Edelsbrunner14}.
Figure~\ref{sfig:Persistence diagram example} shows the persistence diagram of an example data set.
Since the data set, shown in Figure~\ref{sfig:Data set example}, appears to be a ``blob'', the
persistence diagram, as expected, contains few topological features of high persistence in both
dimensions.

\subsection{Intersection Homology and Persistent Intersection Homology}

Despite its prevalence in data analysis, persistent homology exhibits some limitations. In the
context of this paper, we are mostly concerned with its lack of duality for non-manifold data
sets, and with its inability to detect topological features of data sets consisting of
\emph{multiple} manifolds\footnote{%
  We remark that topologically, this case can often be reduced to the computation of ordinary homology,
  because a theorem of Goresky and MacPherson~\cite{Goresky80} ensures that for pseudomanifolds, the
  intersection homology groups remain the same under normalization, and if they are
  nonsingular, the intersection homology groups are ordinary homology groups. As it is not clear how to
  obtain normalizations for real-world data, the calculation of persistent intersection
  homology is necessary.
}.
Recall that for a $d$-manifold, Poincar\'e duality means that the Betti numbers satisfy $\betti{k}
= \betti{d-k}$.
While it is possible to extend persistent homology to obtain something similar for
manifolds~\cite{Cohen-Steiner09,Silva11}, there is no general duality theorem yet.
Additionally, persistent homology cannot detect manifolds of varying dimensionality that are ``glued
together'' in the manner described in Section~\ref{sec:Stratifications}.
For example, we could model the data set from Figure~\ref{fig:Market data}, in which we see a central
``core'' along with some ``flares'', as a topological disk to which we added multiple ``whiskers''.
The persistence diagram does not contain evidence of any whiskers, so the data set will have the same
persistence diagram as a data set that only contains a topological disk. Carlsson~\cite{Carlsson14}
proposes to use filter functions on the data to remedy this situation.
While this helps detect the features, it does not detect that the underlying structure
does not consist of one single manifold.

Intersection homology faces these challenges by providing a homology theory for such spaces with
singularities. We follow the notation of Bendich~\cite{Bendich09,Bendich11} here, who provided
a generic framework for calculating restricted forms of~(persistent) homology, of which intersection
homology is a special case.
In the following, we require a function~$\rest\colon\simplicialcomplex \to \{0,1\}$ that restricts
the usage of simplices. We call a simplex $\sigma$ \emph{proper} or \emph{allowable} if $\rest(\sigma) = 1$.
While $\rest(\simplicialcomplex)$ is not generally a simplicial complex, we can use it
to define a restriction on the chain groups of $\simplicialcomplex$ by calling
a simplicial chain $c \in \chaingroup{p}(\simplicialcomplex)$ \emph{proper} or \emph{allowable} if
both $c$ and $\boundary{p}\,c$ can be written as formal sums of proper simplices.
We refer to the set of allowable $p$-chains as $I^{\rest}\chaingroup{p}(\simplicialcomplex)$.
Since $\boundary{p-1}\circ\boundary{p} = 0$, the boundary of an allowable $p$-chain is an allowable
chain of dimension $p-1$, so the boundary homomorphism gives rise to a chain complex on the set of
allowable chains. We write $I^{\rest}\homologygroup{p}(\simplicialcomplex)$ to denote the $p$\th
homology group of this complex, and refer to it as the $p$\th intersection homology group.
There is a natural restriction of $\rest(\cdot)$ when $\simplicialcomplex$ is filtrated, so we can
define a set of restricted persistent homology groups $I^{\rest}\persistenthomologygroup{p}{i,j}$ in
analogy to the definition of the persistent homology groups.

\subsubsection{Persistent Intersection Homology}
%
To obtain \emph{intersection homology} from this generic framework, we require a few additional
definitions: a \emph{perversity}\footnote{%
  See the unpublished notes by MacPherson on \emph{Intersection Homology and Perverse Sheaves},
  available under \url{http://faculty.tcu.edu/gfriedman/notes/ih.pdf}, for the origin of this name.%
}
is a sequence of integers
\begin{equation}
  \bar{p} = (p_1,p_2,\dots,p_{d-1},p_d)
\end{equation}
such that $-1 \leq p_k \leq k-1$ for every $k$. Alternatively, following the original definition of
Goresky and MacPherson~\cite{Goresky80}, a perversity is a sequence of integers
\begin{equation}
  \bar{p}' = (p_2',p_3',\dots,p_{d-1}',p_d')
\end{equation}
such that $p_2' = 0$ and either $p_{k+1}' = p_k'$ or $p_{k+1}' = p_k' + 1$. Both definitions
permit assessing to what extent a data set deviates from being a manifold.
More precisely, the perversity measures how much deviation from full transverse intersections~(i.e.,
intersections of two submanifolds that yield another submanifold) are permitted for a given simplicial
complex. Each choice of perversity will yield a different set of restricted~(persistent) homology
groups. We focus only on low-dimensional perversities in this paper, with $k \leq 3$.
Finally, tying all the previous definitions together, we define a function~$\rest(\cdot)$ for
a given perversity and a given stratification: a simplex $\sigma$ is considered to be \emph{proper} if
\begin{equation}
  \dim(\sigma \cap X_{d-k}) \leq \dim(\sigma) - k + p_k
  \label{eq:IH}
\end{equation}
holds for all $k \in \{1,\dots,d\}$. Intuitively, this inequality bounds the dimensionality of
the intersection of a simplex with a given subspace. We set $\dim(\emptyset) := -\infty$ so
that simplices without an intersection are considered proper.
Larger values for $p_k$ give us more tolerant intersection conditions, whereas smaller values for
$p_k$ make the intersections more restrictive.
This leads to persistent intersection homology groups with a given perversity function.

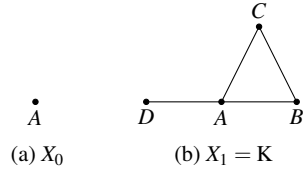
\begin{figure}[t]
  \centering
  \subfloat[$X_0$\label{sfig:Example X0}]{%
    \iffinal
      \includegraphics{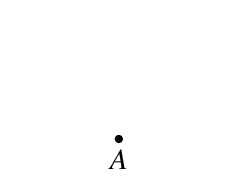}
    \else
      \begin{tikzpicture}
        \coordinate (A) at ( 0.0,0.0);
        \coordinate (B) at ( 1.0,0.0);
        \coordinate (C) at ( 0.5,1.0);
        \coordinate (D) at (-1.0,0.0);

        \node[anchor=north] at (A) {$A$};
        \node[anchor=north, white] at (B) {$B$};
        \node[anchor=south, white] at (C) {$C$};
        \node[anchor=north, white] at (D) {$D$};

        \filldraw (A) circle (1.0pt);
      \end{tikzpicture}
    \fi
  }
  \subfloat[$X_1 = \simplicialcomplex$\label{sfig:Example X1}]{%
    \iffinal
      \includegraphics{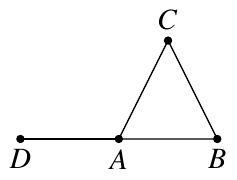}
    \else
      \begin{tikzpicture}
        \coordinate (A) at ( 0.0,0.0);
        \coordinate (B) at ( 1.0,0.0);
        \coordinate (C) at ( 0.5,1.0);
        \coordinate (D) at (-1.0,0.0);

        \node[anchor=north] at (A) {$A$};
        \node[anchor=north] at (B) {$B$};
        \node[anchor=south] at (C) {$C$};
        \node[anchor=north] at (D) {$D$};

        \draw (A) -- (B) -- (C) -- (A);
        \draw (A) -- (D);

        \foreach \c in {A,B,C,D}
          \filldraw (\c) circle (1.0pt);
      \end{tikzpicture}
    \fi
  }
  \caption{%
    A simple example stratified space~\protect\subref{sfig:Example X1} for which simplicial homology is
    incapable of detecting the additional ``whisker''. The \emph{singular
    stratum}~\protect\subref{sfig:Example X0} only consists of a single vertex, $A$.
  }
  \label{fig:Circle with whisker}
\end{figure}

\subsubsection{Simple Example}

Figure~\ref{fig:Circle with whisker} shows a triangulation of a circle with an additional
``whisker''.
This triangulation is in itself not a manifold: at vertex $A$, the neighborhood condition that is
required for a manifold is violated. However, the space is made up of two manifolds, namely a circle
and a line, that are joined at a single point. A natural stratification of such a space thus puts
the singular vertex $A$ in $X_0$ and the full simplicial complex in $X_1$.
With ordinary simplicial homology, we obtain $\betti{0} = 1$, because there is only a single
connected component.
Intersection homology permits only two different perversities here~(we cannot use
Goresky--MacPherson perversities because $d = 1$), either $p_1 = -1$ or $p_1 = 0$; as we are only
interested in $\betti{0}$, we do not have to provide a higher-dimensional value for the perversity.
For $p_1 = -1$, we obtain $\betti{0} = 2$, because \emph{no} simplex that contains $A$ is proper.
This reflects the fact that the simplicial complex is made up of two pieces whose type is different.
For $p_1 = 0$, we obtain again $\betti{0} = 1$ because the singular point now leads to a proper
connected component: Eq.~\ref{eq:IH} becomes  $\dim(\sigma \cap X_1) \leq \dim(\sigma)$, which
is satisfied by every simplex $\sigma$.

\subsubsection{Implementation}

The crucial part of implementing persistent intersection homology lies in an efficient evaluation of
Eq.~\ref{eq:IH}: for each simplex $\sigma$, the calculating the dimension of the
intersection on the left-hand side requires searching through some $X_{d-k}$ and reporting the
intersection with the highest dimension. Large speedups can be obtained by
\begin{inparaenum}[(i)]
  \item restricting the search to $l$-simplices, where $l := \min(\dim\sigma,d-k)$ is the maximum
  dimension that can be achieved by the intersection, and
  \item enumerating all subsets $\tau \subseteq \sigma$~(in reverse lexicographical order, because
  we are looking for the largest dimension) and checking whether $\tau \in X_{d-k}$.
\end{inparaenum}
The second step particularly improves performance when $\dim\sigma$ is small, because we have
to enumerate at most $2^{\dim\sigma}$ simplices and check whether they are part of $X_{d-k}$. Each
check can be done in constant or~(at worst) logarithmic time in the size of $X_{d-k}$.
By contrast, calculating all intersections of $\sigma$ with $X_{d-k}$ takes at least linear time in
the size of $X_{d-k}$. If $2^{\dim\sigma} \cdot \log \card{X_{d-k}} \ll \card{X_{d-k}}$, our method
will be beneficial for performance.
We provide an implementation of persistent intersection homology in
\texttt{Aleph}\footnote{\url{https://github.com/Submanifold/Aleph}}, a software library for
topological data analysis.
We are not aware of any other open-source implementation of persistent
intersection homology at this time.

\section{Using Persistent Intersection Homology}

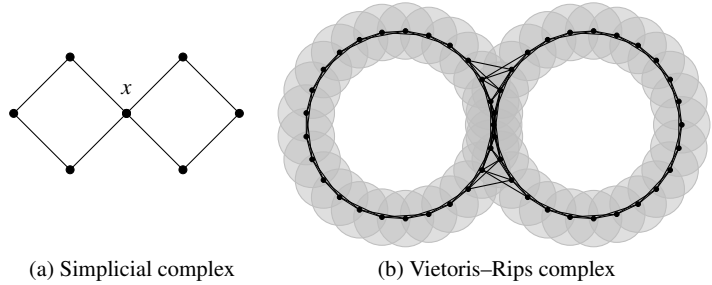
\begin{figure}[tbp]
  \centering
  \subfloat[Simplicial complex\label{sfig:S1vS1 model}]{%
    \raisebox{0.60\height}{%
      \iffinal
        \includegraphics{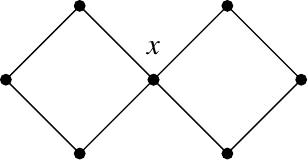}
      \else
        \begin{tikzpicture}[scale=1.50]
          \coordinate (A) at ( 0.0, 0.0);
          \coordinate (B) at ( 0.5, 0.5);
          \coordinate (C) at ( 1.0, 0.0);
          \coordinate (D) at ( 0.5,-0.5);
          \coordinate (E) at (-0.5,-0.5);
          \coordinate (F) at (-1.0, 0.0);
          \coordinate (G) at (-0.5, 0.5);

          \draw (A) -- (B) -- (C) -- (D) -- cycle;
          \draw (A) -- (E) -- (F) -- (G) -- cycle;

          \foreach \c in {A,B,...,G}
            \filldraw (\c) circle (1.0pt);

          \filldraw              (A) circle (1.0pt);
          \node[above=0.15cm] at (A) {$x$};
        \end{tikzpicture}%
      \fi
    }
  }
  \quad
  \subfloat[Vietoris--Rips complex\label{sfig:S1vS1 reality}]{%
    \iffinal
      \includegraphics{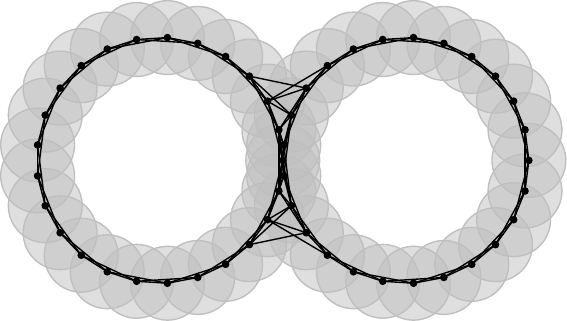}
    \else
      \begin{tikzpicture}[scale=1.25]
        \coordinate (00) at ( 1.0000, 0.0000);
        \coordinate (01) at ( 0.9686, 0.2487);
        \coordinate (02) at ( 0.8763, 0.4818);
        \coordinate (03) at ( 0.7290, 0.6845);
        \coordinate (04) at ( 0.5358, 0.8443);
        \coordinate (05) at ( 0.3090, 0.9511);
        \coordinate (06) at ( 0.0628, 0.9980);
        \coordinate (07) at (-0.1874, 0.9823);
        \coordinate (08) at (-0.4258, 0.9048);
        \coordinate (09) at (-0.6374, 0.7705);
        \coordinate (10) at (-0.8090, 0.5878);
        \coordinate (11) at (-0.9298, 0.3681);
        \coordinate (12) at (-0.9921, 0.1253);
        \coordinate (13) at (-0.9921,-0.1253);
        \coordinate (14) at (-0.9298,-0.3681);
        \coordinate (15) at (-0.8090,-0.5878);
        \coordinate (16) at (-0.6374,-0.7705);
        \coordinate (17) at (-0.4258,-0.9048);
        \coordinate (18) at (-0.1874,-0.9823);
        \coordinate (19) at ( 0.0628,-0.9980);
        \coordinate (20) at ( 0.3090,-0.9511);
        \coordinate (21) at ( 0.5358,-0.8443);
        \coordinate (22) at ( 0.7290,-0.6845);
        \coordinate (23) at ( 0.8763,-0.4818);
        \coordinate (24) at ( 0.9686,-0.2487);
        \coordinate (25) at ( 3.0000, 0.0000);
        \coordinate (26) at ( 2.9686, 0.2487);
        \coordinate (27) at ( 2.8763, 0.4818);
        \coordinate (28) at ( 2.7290, 0.6845);
        \coordinate (29) at ( 2.5358, 0.8443);
        \coordinate (30) at ( 2.3090, 0.9511);
        \coordinate (31) at ( 2.0628, 0.9980);
        \coordinate (32) at ( 1.8126, 0.9823);
        \coordinate (33) at ( 1.5742, 0.9048);
        \coordinate (34) at ( 1.3626, 0.7705);
        \coordinate (35) at ( 1.1910, 0.5878);
        \coordinate (36) at ( 1.0702, 0.3681);
        \coordinate (37) at ( 1.0079, 0.1253);
        \coordinate (38) at ( 1.0079,-0.1253);
        \coordinate (39) at ( 1.0702,-0.3681);
        \coordinate (40) at ( 1.1910,-0.5878);
        \coordinate (41) at ( 1.3626,-0.7705);
        \coordinate (42) at ( 1.5742,-0.9048);
        \coordinate (43) at ( 1.8126,-0.9823);
        \coordinate (44) at ( 2.0628,-0.9980);
        \coordinate (45) at ( 2.3090,-0.9511);
        \coordinate (46) at ( 2.5358,-0.8443);
        \coordinate (47) at ( 2.7290,-0.6845);
        \coordinate (48) at ( 2.8763,-0.4818);
        \coordinate (49) at ( 2.9686,-0.2487);

        \foreach \c in {00,...,49}
          \filldraw[lightgray, fill opacity=0.5] (\twodigits{\c}) circle (0.30cm);

        \draw (00) -- (01); \draw (00) -- (02); \draw (00) -- (23); \draw (00) -- (24);
        \draw (00) -- (36); \draw (00) -- (37); \draw (00) -- (38); \draw (00) -- (39);
        \draw (01) -- (02); \draw (01) -- (03); \draw (01) -- (24); \draw (01) -- (35);
        \draw (01) -- (36); \draw (01) -- (37); \draw (01) -- (38); \draw (02) -- (03);
        \draw (02) -- (04); \draw (02) -- (34); \draw (02) -- (35); \draw (02) -- (36);
        \draw (02) -- (37); \draw (03) -- (04); \draw (03) -- (05); \draw (03) -- (35);
        \draw (03) -- (36); \draw (04) -- (05); \draw (04) -- (06); \draw (05) -- (06);
        \draw (05) -- (07); \draw (06) -- (07); \draw (06) -- (08); \draw (07) -- (08);
        \draw (07) -- (09); \draw (08) -- (09); \draw (08) -- (10); \draw (09) -- (10);
        \draw (09) -- (11); \draw (10) -- (11); \draw (10) -- (12); \draw (11) -- (12);
        \draw (11) -- (13); \draw (12) -- (13); \draw (12) -- (14); \draw (13) -- (14);
        \draw (13) -- (15); \draw (14) -- (15); \draw (14) -- (16); \draw (15) -- (16);
        \draw (15) -- (17); \draw (16) -- (17); \draw (16) -- (18); \draw (17) -- (18);
        \draw (17) -- (19); \draw (18) -- (19); \draw (18) -- (20); \draw (19) -- (20);
        \draw (19) -- (21); \draw (20) -- (21); \draw (20) -- (22); \draw (21) -- (22);
        \draw (21) -- (23); \draw (22) -- (23); \draw (22) -- (24); \draw (22) -- (39);
        \draw (22) -- (40); \draw (23) -- (24); \draw (23) -- (38); \draw (23) -- (39);
        \draw (23) -- (40); \draw (23) -- (41); \draw (24) -- (37); \draw (24) -- (38);
        \draw (24) -- (39); \draw (24) -- (40); \draw (25) -- (26); \draw (25) -- (27);
        \draw (25) -- (48); \draw (25) -- (49); \draw (26) -- (27); \draw (26) -- (28);
        \draw (26) -- (49); \draw (27) -- (28); \draw (27) -- (29); \draw (28) -- (29);
        \draw (28) -- (30); \draw (29) -- (30); \draw (29) -- (31); \draw (30) -- (31);
        \draw (30) -- (32); \draw (31) -- (32); \draw (31) -- (33); \draw (32) -- (33);
        \draw (32) -- (34); \draw (33) -- (34); \draw (33) -- (35); \draw (34) -- (35);
        \draw (34) -- (36); \draw (35) -- (36); \draw (35) -- (37); \draw (36) -- (37);
        \draw (36) -- (38); \draw (37) -- (38); \draw (37) -- (39); \draw (38) -- (39);
        \draw (38) -- (40); \draw (39) -- (40); \draw (39) -- (41); \draw (40) -- (41);
        \draw (40) -- (42); \draw (41) -- (42); \draw (41) -- (43); \draw (42) -- (43);
        \draw (42) -- (44); \draw (43) -- (44); \draw (43) -- (45); \draw (44) -- (45);
        \draw (44) -- (46); \draw (45) -- (46); \draw (45) -- (47); \draw (46) -- (47);
        \draw (46) -- (48); \draw (47) -- (48); \draw (47) -- (49); \draw (48) -- (49);

        \foreach \c in {00,...,49}
          \filldraw (\twodigits{\c}) circle (0.025cm);
      \end{tikzpicture}
    \fi
  }
  \caption{%
    Calculating the Vietoris--Rips complex of a point cloud makes it impossible to detect
    singularities by homological means alone.
  }
  \label{fig:S1vS1}
\end{figure}

Prior to using persistent intersection homology in a topological data analysis workflow, we need to
discuss one of its pitfalls: the Vietoris--Rips complex is commonly used in topological data
analysis to deal with multivariate data sets.
For persistent intersection homology, this construction turns out to result in triangulations that
yield unexpected results.
Figure~\ref{fig:S1vS1} depicts an example of this issue. Here we see the \emph{one-point
union}, i.e., the wedge sum, of two circles, denoted by $S^1 \vee S^1$.
Formally, this can be easily modeled as a simplicial complex~$\simplicialcomplex$~(Figure~\ref{sfig:S1vS1 model}).
The smallest stratification of this space places the singular point~$x$ in its own subspace, i.e.,
$X_0 = \{x\}$, $X_1 = \simplicialcomplex$, and uses $\bar{p} = (-1)$.
The intersection homology of $\simplicialcomplex$  results in $\betti{0} = 2$, because of the singular
point at which the two circles are connected. Calculating persistent intersection homology of
a point cloud that describes this space~(Figure~\ref{sfig:S1vS1 reality}), by contrast, results in
$\betti{0} = 1$, regardless of whether we ensure that the triangulation is
\emph{flaglike}~\cite{MacPherson86} by performing the first barycentric subdivision~(which is
guaranteed to make the calculations \emph{independent} of the stratification~\cite{Goresky80}).
The reason for this is that the topological realization of the Vietoris--Rips complex seems to be
more closely tied to regular neighborhoods than to the homeomorphism type of $S^1 \vee S^1$.
However, the regular neighborhood of a space is always a manifold. It can be thought of as
calculating a ``thickened'' version of the space in which isolated singularities disappear.

As far as we know, Bendich and Harer~\cite{Bendich11}, while discussing other dependencies of
persistent intersection homology, did not discuss this aspect. Yet, it is crucial to get persistent
intersection homology to ``detect'' those singularities if we want to understand the manifold
structure of a given data set. To circumvent this issue, we propose obtaining additional information
about the geometry of a given point cloud in order to determine which points are supposed to be
singular. Alternatively, we could try to learn a suitable stratification of the whole
space~\cite{Bendich12} at the cost of reduced performance.

\subsection{Choosing a Stratification}

Having seen that the utility and expressiveness of persistent intersection homology hinge upon the
choice of a stratification, we now develop several constructions.
We restrict ourselves to the detection of isolated singular points, i.e., vertices or $0$-simplices,
in this paper. A stratification should ideally reflect the existence of singularities in a data set.
For the example shown in Figure~\ref{fig:Circle with whisker}, a singularity exists at $A$ because
the ``whisker'' will remain a one-dimensional piece regardless of the scale at which we look at the
data, while the triangle is a two-dimensional object.
This observation leads to a set of stratification strategies, which we first detail before applying
them in Section~\ref{sec:Results}.

\subsubsection{Dimensionality-Based Stratifications}
%
In order to stratify unstructured data according to the local intrinsic dimensionality, we propose
the following scheme. We first obtain the $k$ nearest neighbors of every data point and treat them
as local patches.
For each of these subsets, we perform a principal component analysis~(PCA) and obtain the respective set of eigenvalues
$\{ \lambda_1, \dots, \lambda_d \}$, where $d$ refers to the maximum number of attributes in the
point cloud.
We then calculate the largest spectral gap, i.e.,
\begin{equation}
  d_i := \argmax_{j \in \{2,\dots,d\}} | \lambda_j - \lambda_{j-1} | - 1,
\end{equation}
and use it as an estimate of the local intrinsic dimensionality at the $i$\th data point.
Points that can be well represented by a single eigenvalue are thus taken to correspond to a locally
one-dimensional patch in the data, for example.
In practice, as PCA is not robust against outliers, one typically requires some smoothing iterations
for the estimates. We use several iterations of smoothing based on nearest neighbors, similar to
\emph{mean shift clustering}~\cite{Cheng95}.
The resulting values can then be used to stratify according to local dimensionality.

\subsubsection{Density-based Stratifications}
%
We can also stratify unstructured data according to the behavior of a density estimator, such as
a \emph{truncated Gaussian kernel}, i.e.,
\begin{equation}
  f(x) := \mathlarger{\mathlarger\sum}_{y \neq x} \exp{\bigg( - \frac{ \dist^2(x,y) }{ 2 h }\bigg)},
\end{equation}
where $h$ is the bandwidth of the estimator and we define the exponential expression to be $0$ if $\dist(x,y) > h$.
The density values give rise to a distribution of values so that we can use standard outlier
detection methods. Once outliers have been identified, they can be put into the first subset of the
filtration. This approach has the advantage of rapidly detecting interesting data points but it
cannot be readily extended to higher-dimensional simplices.

\subsubsection{Curvature-Based Stratifications}

The curvature of a manifold is an important property that can be used to detect differences in local
structure. Using a standard algorithm to estimate curvature in meshes~\cite{Meyer03}, we can easily
identify a region around the singular point in the ``pinched torus'' as having an extremely small
curvature. Figure~\ref{fig:Pinched torus curvature} depicts this. For higher-dimensional point
clouds, we propose obtaining an approximation of curvature by using the curvature of
high-dimensional spheres that are fit to local patches of a point cloud.
More precisely, we extract the $k$ nearest neighbors of every point in a point cloud and fit
a high-dimensional sphere. Such a fit can be accomplished using standard least squares approaches,
such as the one introduced by Pratt~\cite{Pratt87}.

\section{Results}
\label{sec:Results}

In the following, we discuss the benefits of persistent intersection homology over ordinary
persistent homology by means of several data sets, containing random samples of non-trivial
topological pseudomanifolds, as well as experimental data from image processing.

\subsection{Wedge of Spheres}

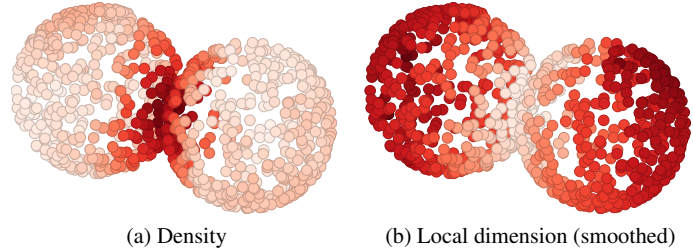
\begin{figure}[tbp]
  \centering
  \subfloat[Density]{%
    \iffinal
      \includegraphics{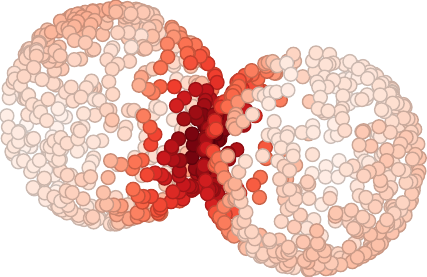}
    \else
      \begin{tikzpicture}
        \begin{axis}[%
          axis lines        = none,
          grid              = none,
          z buffer          = sort,
          unit vector ratio*= 1 1 1,
          colormap/Reds,
          point meta        = explicit,
        ]
          \addplot3[only marks, scatter] file {Data/S2vS2/Example_point_cloud_with_density.txt};
        \end{axis}
      \end{tikzpicture}
    \fi
  }%
  \quad
  \subfloat[Local dimension~(smoothed)]{%
    \iffinal
      \includegraphics{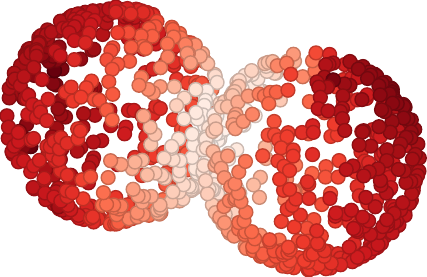}
    \else
      \begin{tikzpicture}
        \begin{axis}[%
          axis lines        = none,
          grid              = none,
          z buffer          = sort,
          unit vector ratio*= 1 1 1,
          colormap/Reds,
          point meta        = explicit,
        ]
          \addplot3[only marks, scatter] file {Data/S2vS2/Example_point_cloud_with_local_dimensionality.txt};
        \end{axis}
      \end{tikzpicture}
    \fi
  }
  \caption{%
    A random sample of $S^2 \vee S^2$, color-coded by two stratification strategies.  Both
    descriptors register either extremely high~(density) or extremely low~(dimensionality) values as
    we approach the singular part of the data set. The corresponding points are put into $X_0$.
  }
  \label{fig:S2vS2}
\end{figure}

\begin{figure}[tbp]
  \centering
  \subfloat[Persistent homology\label{sfig:S2vS2 barcode PH}]{%
    \iffinal
      \includegraphics{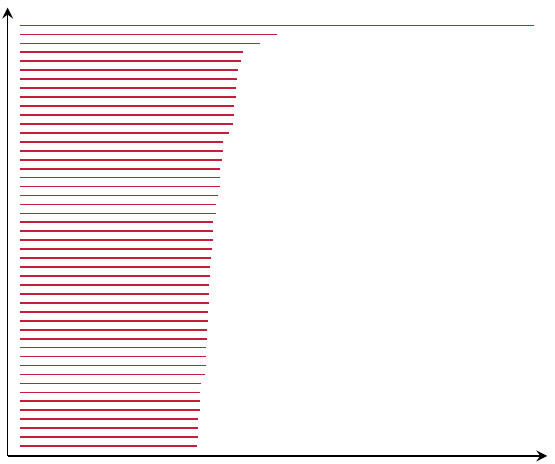}
    \else
      \begin{tikzpicture}
        \begin{axis}[
          axis x line         = bottom,
          axis y line         = left,
          enlarge y limits    = {abs=0.100, lower},
          enlarge x limits    = 0.025,
          empty line          = jump,
          ymin                = 800,
          ymax                = 850,
          restrict y to domain= 801:850,
          yticklabels         = {,,},
          xticklabels         = {,,},
          ticks               = none,
          scale               = 0.80,
        ]
          \addplot[cardinal, no marks,line width=0.5pt] table {./Data/S2vS2/PH_0_barcode.txt};
        \end{axis}
        \path ([shift={(-5\pgflinewidth,-5\pgflinewidth)}]current bounding box.south west)
              ([shift={( 5\pgflinewidth, 5\pgflinewidth)}]current bounding box.north east);
      \end{tikzpicture}%
    \fi
  }
  \subfloat[Persistent intersection homology\label{sfig:S2vS2 barcode IH}]{%
    \iffinal
      \includegraphics{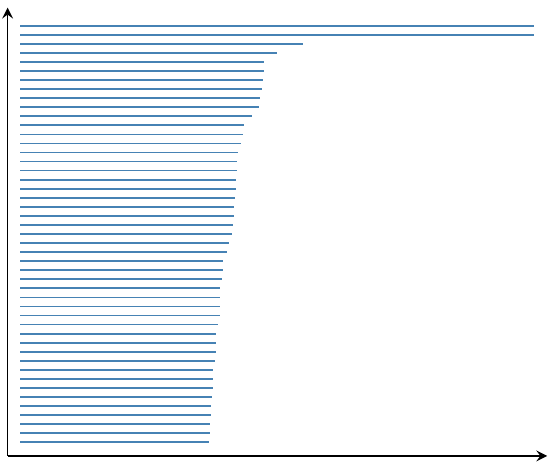}
    \else
      \begin{tikzpicture}
        \begin{axis}[
          axis x line          = bottom,
          axis y line          = left,
          enlarge y limits     = {rel=0.100, lower},
          enlarge x limits     = 0.025,
          empty line           = jump,
          ymin                 = 750,
          ymax                 = 795,
          restrict y to domain = 747:795,
          yticklabels          = {,,},
          xticklabels          = {,,},
          ticks                = none,
          scale                = 0.80,
        ]
          \addplot[yale,mark=none,line width=0.5pt] table {./Data/S2vS2/IH_0_barcode.txt};
        \end{axis}
        \path ([shift={(-5\pgflinewidth,-5\pgflinewidth)}]current bounding box.south west)
              ([shift={( 5\pgflinewidth, 5\pgflinewidth)}]current bounding box.north east);
      \end{tikzpicture}
    \fi
  }
  \caption{%
    Excerpt of the zero-dimensional barcodes for $S^2 \vee S^2$.
    With persistent homology~\protect\subref{sfig:S2vS2 barcode PH}, no additional connected
    component appears, whereas with persistent intersection homology~\protect\subref{sfig:S2vS2
    barcode IH} with the density-based stratification, the singular point/region results in
    splitting the data.
  }
  \label{fig:S2vS2 barcodes}
\end{figure}
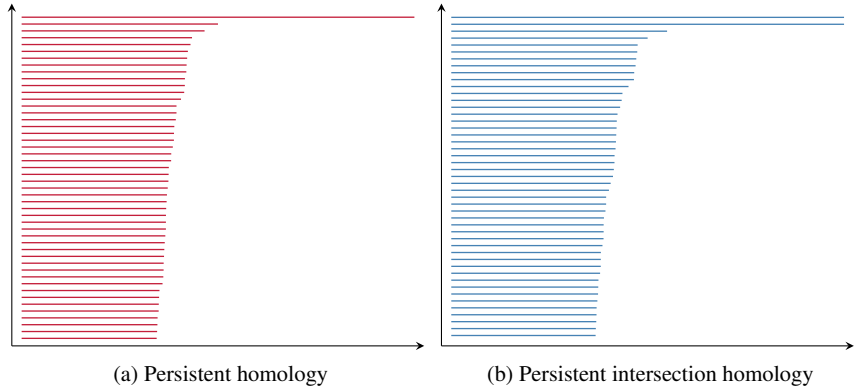

We extend the example depicted in Figure~\ref{fig:S1vS1} and sample points at random from a wedge of
$2$-spheres. If no precautions are taken, the resulting data set suffers from the problem that we previously
outlined. We thus use it to demonstrate the efficacy of our stratification strategies.
Figure~\ref{fig:S2vS2} depicts the data set along with two different descriptors.
In both cases, we build a simple stratification in which $X_0$ contains all singular points, $X_1
= X_0$, and $X_2 = \simplicialcomplex$, i.e., the original space. We use the default
Goresky--MacPherson perversity $\bar{p}' = (0)$. This suffices to detect that the data set is not
a manifold: we obtain $\betti{0} = 2$ for both stratification strategies, whereas persistent
homology only shows $\betti{0} = 1$. Figure~\ref{fig:S2vS2 barcodes} depicts excerpts of the
zero-dimensional barcodes for the data set.
The two topological features with infinite persistence are clearly visible in the persistent
intersection homology barcode.
Since $\betti{2} = 2$, this re-establishes Poincar\'e duality.

\subsection{Pinched Torus}

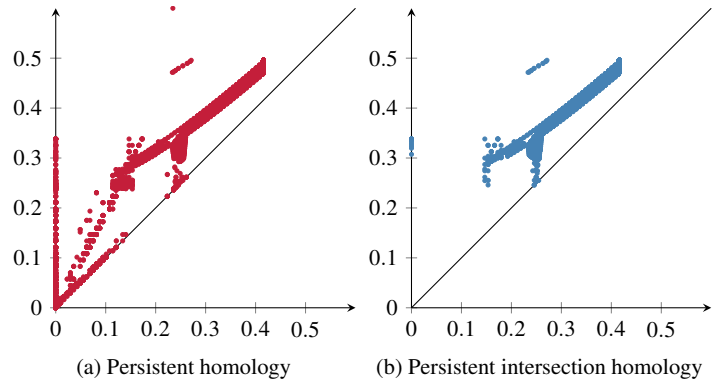
\begin{figure}[tbp]
  \centering
  \subfloat[Persistent homology\label{sfig:Pinched torus PH}]{%
    \iffinal
      \includegraphics{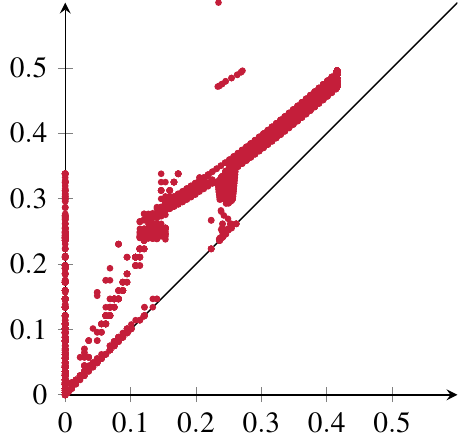}
    \else
      \begin{tikzpicture}
        \begin{axis}[%
          axis x line       = bottom,
          axis y line       = left,
          enlargelimits     = false,
          mark size         = 0.75pt,
          unit vector ratio*= 1 1 1,
          enlargelimits     = false,
          xmin              = 0,
          ymin              = 0,
          scale             = 0.70,
          xtick             = {0, 0.1, 0.2, 0.3, 0.4, 0.5},
          ytick             = {0, 0.1, 0.2, 0.3, 0.4, 0.5},
        ]
          \addplot[cardinal, only marks] file {Data/Pinched_torus/PH_1.txt};
          \addplot[domain=0.0:0.6]            {x};
        \end{axis}
      \end{tikzpicture}
    \fi
  }
  \subfloat[Persistent intersection homology\label{sfig:Pinched torus IH}]{%
    \iffinal
      \includegraphics{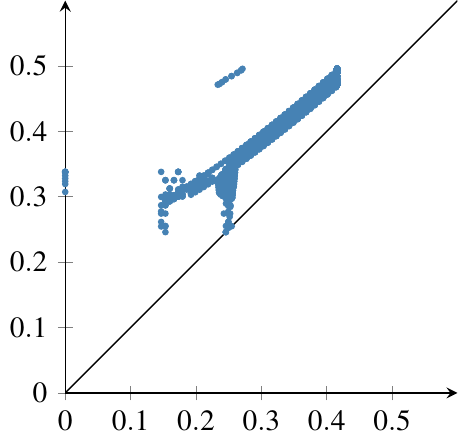}
    \else
      \begin{tikzpicture}
        \begin{axis}[%
          axis x line       = bottom,
          axis y line       = left,
          enlargelimits     = false,
          mark size         = 0.75pt,
          unit vector ratio*= 1 1 1,
          enlargelimits     = false,
          xmin              = 0,
          ymin              = 0,
          scale             = 0.70,
          xtick             = {0, 0.1, 0.2, 0.3, 0.4, 0.5},
          ytick             = {0, 0.1, 0.2, 0.3, 0.4, 0.5},
        ]
          \addplot[yale, only marks] file {Data/Pinched_torus/IH_1.txt};
          \addplot[domain=0.0:0.6]            {x};
        \end{axis}
      \end{tikzpicture}
    \fi
  }
  \caption{%
    Persistent homology detects more one-dimensional features for the ``pinched torus'' data set
    than persistent intersection homology.
  }
  \label{fig:Pinched torus persistence diagrams}
\end{figure}

We demonstrate the curvature-based stratification using the ``pinched
torus'' data set.
Figure~\ref{fig:Pinched torus curvature} depicts the torus along with curvature estimates.
A standard outlier test helps us detect the region around the singular point. We set up the
stratification such that $X_0$ contains all points from the detected region, $X_1 = X_0$, and $X_2
= \simplicialcomplex$. Moreover, we use $\bar{p}' = (0)$ because the dimensionality of the input data
prevents us from detecting any higher-dimensional features.
Persistent homology shows that the point cloud contains a persistent cycle in dimension one.
Essentially, the data are considered to be a ``thickened circle''.
Figure~\ref{fig:Pinched torus persistence diagrams} depicts the persistence diagrams. We can see
that the point with infinite persistence~(shown in Figure~\ref{sfig:Pinched torus PH} at the top
border) is \emph{missing} in addition to many other points in the persistent intersection homology
diagram~(Figure~\ref{sfig:Pinched torus IH}). The Wasserstein distance~\cite{Edelsbrunner10} between the two diagrams is
thus large, indicating the non-manifold structure of the data.

\subsection{Synthetic Faces}

This data set was originally used to demonstrate the effectiveness of nonlinear dimensionality
reduction algorithms~\cite{Tenenbaum00}. Previous research demonstrated that the data set does not
exhibit uniform density~\cite{Rieck15b}, which makes the existence of~(isolated) singular points
possible. It is known that the intrinsic dimension of the data set is three, so we shall only take
a look at low-dimensional topological features.
More precisely, using the curvature-based stratification, we want to see how persistence diagrams in
dimensions 0--2 change when we calculate intersection homology.
Note that analyzing three-dimensional features is not expedient, because the stratification cannot
detect deviations from ``manifoldness'' in this dimension.

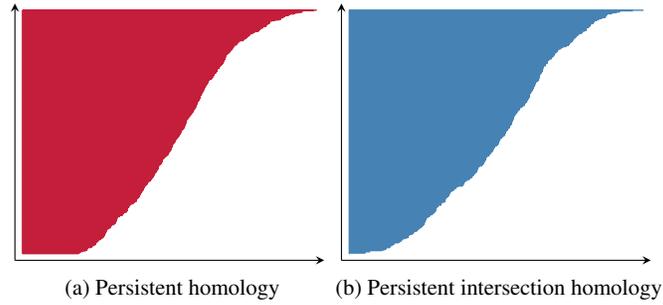
\begin{figure}[tbp]
  \centering
  \subfloat[Persistent homology\label{sfig:Faces barcode PH}]{%
    \iffinal
      \includegraphics{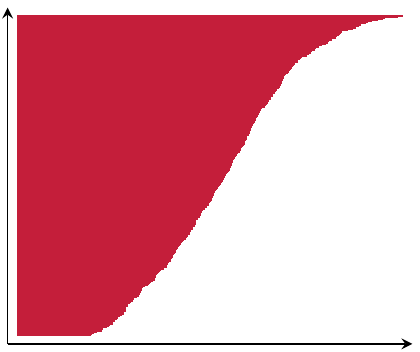}
    \else
      \begin{tikzpicture}
        \begin{axis}[
          axis x line          = bottom,
          axis y line          = left,
          enlarge y limits     = 0.025,
          enlarge x limits     = 0.025,
          empty line           = jump,
          ymin                 = 150,
          restrict y to domain = 151:600,
          yticklabels          = {,,},
          xticklabels          = {,,},
          ticks                = none,
          scale                = 0.60,
        ]
          \addplot[cardinal, no marks,line width=0.5pt] table {./Data/Faces/MC/PH_0_barcode.txt};
        \end{axis}
        \path ([shift={(-5\pgflinewidth,-5\pgflinewidth)}]current bounding box.south west)
              ([shift={( 5\pgflinewidth, 5\pgflinewidth)}]current bounding box.north east);
      \end{tikzpicture}
    \fi
  }
  \subfloat[Persistent intersection homology\label{sfig:Faces barcode IH}]{%
    \iffinal
      \includegraphics{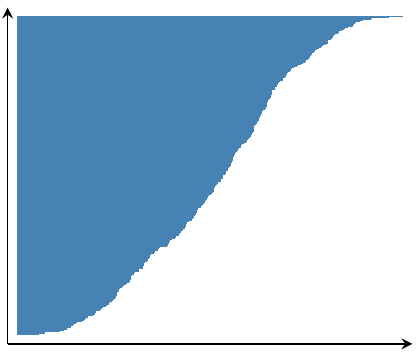}
    \else
      \begin{tikzpicture}
        \begin{axis}[
          axis x line          = bottom,
          axis y line          = left,
          enlarge y limits     = 0.025,
          enlarge x limits     = 0.025,
          empty line           = jump,
          ymin                 = 100,
          restrict y to domain = 101:560,
          yticklabels          = {,,},
          xticklabels          = {,,},
          ticks                = none,
          scale                = 0.60,
        ]
          \addplot[yale,mark=none,line width=0.5pt] table {./Data/Faces/MC/IH_0_barcode.txt};
        \end{axis}
        \path ([shift={(-5\pgflinewidth,-5\pgflinewidth)}]current bounding box.south west)
              ([shift={( 5\pgflinewidth, 5\pgflinewidth)}]current bounding box.north east);
      \end{tikzpicture}
    \fi
  }
  \caption{%
    Zero-dimensional barcodes for the ``Synthetic Faces'' data set. Both barcodes are virtually
    identical, indicating that the singular points do not influence connected components.
  }
  \label{fig:Faces barcodes}
\end{figure}

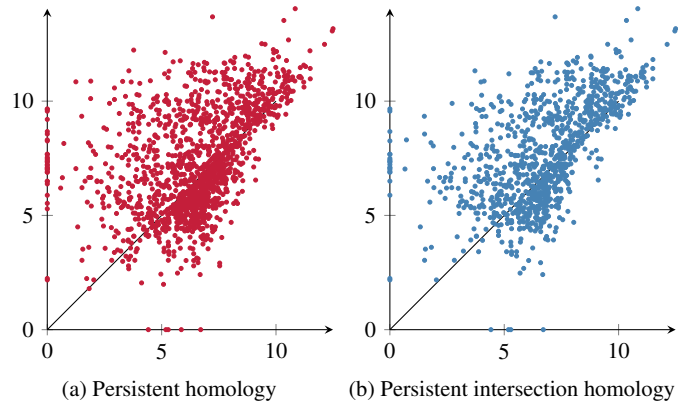
\begin{figure}[tbp]
  \centering
  \subfloat[Persistent homology]{
    \iffinal
      \includegraphics{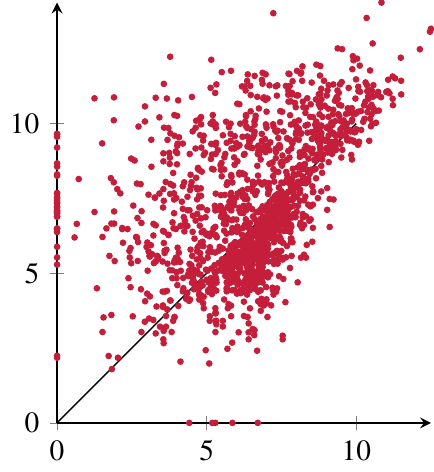}
    \else
      \begin{tikzpicture}
        \begin{axis}[%
          axis x line       = bottom,
          axis y line       = left,
          enlargelimits     = false,
          mark size         = 0.75pt,
          unit vector ratio*= 1 1 1,
          enlargelimits     = false,
          xmin              = 0,
          ymin              = 0,
          scale             = 0.75,
        ]
          \addplot[cardinal, only marks] file                        {Data/Faces/MC/PH_1.txt};
          \addplot[cardinal, only marks] table[x index=1, y index=0] {Data/Faces/MLE/PH_2.txt};
          \addplot[domain=0.0:10]                                    {x};
        \end{axis}
      \end{tikzpicture}
    \fi
  }
  \subfloat[Persistent intersection homology]{
    \iffinal
      \includegraphics{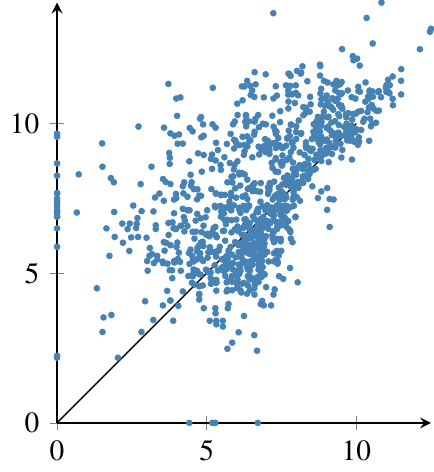}
    \else
      \begin{tikzpicture}
        \begin{axis}[%
          axis x line       = bottom,
          axis y line       = left,
          enlargelimits     = false,
          mark size         = 0.75pt,
          unit vector ratio*= 1 1 1,
          enlargelimits     = false,
          xmin              = 0,
          ymin              = 0,
          scale             = 0.75,
        ]
          \addplot[yale, only marks] file                        {Data/Faces/MC/IH_1.txt};
          \addplot[yale, only marks] table[x index=1, y index=0] {Data/Faces/MLE/IH_2.txt};
          \addplot[domain=0.0:10]                                {x};
        \end{axis}
      \end{tikzpicture}
    \fi
  }
  \caption{%
    Comparison of persistent homology in dimension one~(above diagonal) and two~(below diagonal) for
    the ``Synthetic Faces'' data set. The overall structure is similar, and only few features
    disappear during the calculation of persistent intersection homology.
  }
  \label{fig:Faces dimension 1}
\end{figure}

Figure~\ref{fig:Faces barcodes} depicts the zero-dimensional barcodes of the data set. They are
virtually identical for both methods~(we find that their Wasserstein distance is extremely small),
except for some minor shifts in the destruction values, i.e., the endpoints of every interval.
This indicates that the singular points only have a very \emph{local} influence on the structure of
the data set; they are not resulting in a split, for example.
For dimensions one and two, depicted by Figure~\ref{fig:Faces dimension 1}, we observe a similar behavior.
The overall structure of both persistence diagrams is similar, and there is only a slight decrease
in total persistence~\cite{Cohen-Steiner10} for persistent intersection homology.
Likewise, the Wasserstein distance between both diagrams is extremely small.

In summary, we see that we are unable to detect significant differences in zero-dimensional,
one-dimensional, and two-dimensional topological features. This lends credibility to the assumption
that the data set is a \emph{single} manifold.

\section{Conclusion}

We showed how to use persistent intersection homology for the analysis of data sets that might not
represent a single manifold. Moreover, we described some pitfalls when applying this
technique---namely, finding suitable stratifications, and presented several strategies for doing so.
We demonstrated the utility of persistent intersection homology on several data sets of
low intrinsic dimensionality.
Future work could focus on improving the performance of the admissibility condition in
Eq.~\ref{eq:IH} to process data sets with higher intrinsic dimensions.
It would also be interesting to extend stratification strategies to higher-dimensional
strata, i.e., singular \emph{regions} instead of singular \emph{points}.


\bibliographystyle{spmpsci}
\bibliography{TopoInVis2017_Intersection_Homology}

\begin{thebibliography}{10}
\providecommand{\url}[1]{{#1}}
\providecommand{\urlprefix}{URL }
\expandafter\ifx\csname urlstyle\endcsname\relax
  \providecommand{\doi}[1]{DOI~\discretionary{}{}{}#1}\else
  \providecommand{\doi}{DOI~\discretionary{}{}{}\begingroup
  \urlstyle{rm}\Url}\fi

\bibitem{Bendich09}
Bendich, P.: Analyzing stratified spaces using persistent versions of
  intersection and local homology.
\newblock Ph.D. thesis, Duke University (2009)

\bibitem{Bendich11}
Bendich, P., Harer, J.: Persistent intersection homology.
\newblock FoCM \textbf{11}(3), 305--336 (2011)

\bibitem{Bendich12}
Bendich, P., Wang, B., Mukherjee, S.: Local homology transfer and
  stratification learning.
\newblock In: Y.~Rabani (ed.) Symposium on Discrete Algorithms, pp. 1355--1370.
  SIAM (2012)

\bibitem{Carlsson14}
Carlsson, G.: Topological pattern recognition for point cloud data.
\newblock Acta Numerica \textbf{23}, 289--368 (2014)

\bibitem{Cheng95}
Cheng, Y.: Mean shift, mode seeking, and clustering.
\newblock IEEE TPAMI \textbf{17}(8), 790--799 (1995)

\bibitem{Cohen-Steiner07}
Cohen-Steiner, D., Edelsbrunner, H., Harer, J.: Stability of persistence
  diagrams.
\newblock Discrete {\&} Computational Geometry \textbf{37}(1), 103--120 (2007)

\bibitem{Cohen-Steiner09}
Cohen-Steiner, D., Edelsbrunner, H., Harer, J.: Extending persistence using
  {P}oincar{\'e} and {L}efschetz duality.
\newblock FoCM \textbf{9}(1), 79--103 (2009)

\bibitem{Cohen-Steiner10}
Cohen-Steiner, D., Edelsbrunner, H., Harer, J., Mileyko, Y.: Lipschitz
  functions have {$\mathrm{L}_p$}-stable persistence.
\newblock Foundations of Computational Mathematics \textbf{10}(2), 127--139
  (2010)

\bibitem{Donoho05}
Donoho, D.L., Grimes, C.: Image manifolds which are isometric to {E}uclidean
  space.
\newblock Journal of Mathematical Imaging and Vision \textbf{23}(1), 5--24
  (2005)

\bibitem{Edelsbrunner10}
Edelsbrunner, H., Harer, J.: Computational topology: {A}n introduction.
\newblock AMS (2010)

\bibitem{Edelsbrunner14}
Edelsbrunner, H., Morozov, D.: Persistent homology: {T}heory and practice.
\newblock In: European Congress of Mathematics. EMS Publishing House, Z\"urich,
  Switzerland (2014)

\bibitem{Fefferman16}
Fefferman, C., Mitter, S., Narayanan, H.: Testing the manifold hypothesis.
\newblock Journal of the American Mathematical Society \textbf{29}(4),
  983--1049 (2016)

\bibitem{Goresky80}
Goresky, M., MacPherson, R.: Intersection homology theory.
\newblock Topology \textbf{19}(2), 135--162 (1980)

\bibitem{Hinton97}
Hinton, G.E., Dayan, P., Revow, M.: Modeling the manifolds of images of
  handwritten digits.
\newblock IEEE Transactions on Neural Networks \textbf{8}(1), 65--74 (1997)

\bibitem{Kirwan06}
Kirwan, F., Woolf, J.: An Introduction to Intersection Homology Theory, 2{\nd}
  edn.
\newblock Chapman and Hall/CRC, Boca Raton, FL, USA (2006)

\bibitem{MacPherson86}
MacPherson, R., Vilonen, K.: Elementary construction of perverse sheaves.
\newblock Inventiones mathematicae \textbf{84}(2), 403--435 (1986)

\bibitem{Meyer03}
Meyer, M., Desbrun, M., Schr{\"o}der, P., Barr, A.H.: Discrete
  differential-geometry operators for triangulated {$2$}-manifolds.
\newblock In: H.C. Hege, K.~Polthier (eds.) Visualization and Mathematics III,
  pp. 35--57. Springer, Heidelberg (2003)

\bibitem{Narayanan10}
Narayanan, H., Mitter, S.: Sample complexity of testing the manifold
  hypothesis.
\newblock In: NIPS~23, pp. 1786--1794. Curran Associates, Inc., Red Hook, NY,
  USA (2010)

\bibitem{Pratt87}
Pratt, V.: Direct least-squares fitting of algebraic surfaces.
\newblock ACM SIGGRAPH Computer Graphics \textbf{21}(4), 145--152 (1987)

\bibitem{Rieck15b}
Rieck, B., Leitte, H.: Persistent homology for the evaluation of dimensionality
  reduction schemes.
\newblock Computer Graphics Forum \textbf{34}(3), 431--440 (2015)

\bibitem{Roweis00}
Roweis, S.T., Saul, L.K.: Nonlinear dimensionality reduction by locally linear
  embedding.
\newblock Science \textbf{290}(5500), 2323--2326 (2000)

\bibitem{Saul03}
Saul, L.K., Roweis, S.T.: Think globally, fit locally: {U}nsupervised learning
  of low dimensional manifolds.
\newblock Journal of Machine Learning Research \textbf{4}, 119--155 (2003)

\bibitem{Silva11}
de~Silva, V., Morozov, D., Vejdemo-Johansson, M.: Dualities in persistent
  (co)homology.
\newblock Inverse Problems \textbf{27}(12), 124,003 (2011)

\bibitem{Singh07}
Singh, G., M{\'e}moli, F., Carlsson, G.: Topological methods for the analysis
  of high dimensional data sets and {3D} object recognition.
\newblock In: Eurographics Symposium on Point-Based Graphics. Eurographics
  Association (2007)

\bibitem{Tenenbaum00}
Tenenbaum, J.B., {de Silva}, V., Langford, J.C.: A global geometric framework
  for nonlinear dimensionality reduction.
\newblock Science \textbf{290}(5500), 2319--2323 (2000)

\end{thebibliography}

\end{document}